# The Fixed Point Property
# of Quasi-Point-Separable Topological Vector Spaces

Jinlu Li

Department of Mathematics
Shawnee State University
Portsmouth, Ohio 45662
USA

**Dedicated to Prof. Dr. Christiane Tammer in the occasion of her 65th birthday**

## Abstract

In this paper, we introduce the concept of quasi-point-separable topological vector spaces, which has the following important properties:

1. In general, the conditions for a topological vector space to be quasi-point-separable is not very difficult to check;
2. The class of quasi-point-separable topological vector spaces is very large that includes locally convex topological vector spaces and pseudonorm adjoint topological vector spaces as special cases;
3. Every quasi-point-separable Housdorrf topological vector space has the fixed point property (that is, every continuous self-mapping on any given nonempty closed and convex subset has a fixed point), which is the result of the main theorem of this paper (Theorem 4.1);

Furthermore, we provide some concrete examples of quasi-point-separable topological vector spaces, which are not locally convex. It follows that the main theorem of this paper is a proper extension of Tychonoff's fixed point theorem on locally convex topological vector spaces.

1. Introduction

In fixed point theory, a topological space is said to have the fixed point property if every continuous self-map on this space has at least one fixed point (see [7]).

In this paper, we say that a topological vector space $(X, \tau)$ has the fixed point property if every non-empty, compact and convex subset of $X$ has the fixed point property.

Since the fixed point property is topological invariant, i.e. it is preserved by any homeomorphism. And the fixed point property is also preserved by any retraction. Then, according to Brouwer fixed point theorem, we have

**An alternative version of Brouwer's Fixed Point Theorem** (see [1]). *Every n-dimension Euclidean space $\mathbb{R}^n$ has the fixed point property.*

In 1930, Schauder extended Brouwer's fixed-point theorem from Euclidean spaces to Banach spaces.

**Schauder's fixed-point theorem (an alternative version [15]).** *Every Banach space has the fixed point property.*

To generalize the underlying spaces for fixed point property, in 1934, Tychonoff extended Schauder's fixed point theorem from Banach spaces to locally convex topological vector spaces.

**Tychonoff's fixed point theorem (an alternative version [18]).** *Every Hausdorff locally convex topological vector space has the fixed point property.*

After Schauder proved his fixed point theorem, Schauder raised a well-known conjecture in the field of fixed point theory.

**Schauder's Conjecture (an alternative version)**. *Every Hausdorff topological vector space has the fixed point property.*

Since then, the Schauder's Conjecture has become one of the most important open problems in the field of nonlinear analysis. The Schauder's Conjecture and related fixed point problems have attracted many authors' attention (see [2−6], [11−14], [16−17]). In 2001, R. Cauty in [3] proposed an answer to the Schauder's Conjecture. After Cauty's paper published for four years, in the International Conference of Fixed Point Theory and its Applications in 2005, Professor T. Dobrowolski remarked that there is a gap in Cauty's proof. Therefore, Schauder's Conjecture is still remain unsolved.

In 2020, the present author introduced the concept of pseudonorm adjoint topological vector spaces; and proved that every Housdorrf and total pseudonorm adjoint topological vector space has the fixed point property. It is stated as below.

**An alternative version of a fixed point theorem in [9] (2020).** *Suppose that $(X, \tau)$ is a Housdorrf and total pseudonorm adjoint topological vector space. Then X has the fixed point property.*

Meanwhile a concrete example was provided in this paper [9] showing that the concept of pseudonorm adjoint topological vector spaces is a proper extension of locally convex topological vector spaces; and therefore, the above fixed point theorem is a proper extension of the Tychonoff's fixed point theorem.

In this paper, we introduce a new concept of quasi-point-separable topological vector spaces, which is a generalization of pseudonorm adjoint topological vector spaces. It immediately implies that is a proper extension of locally convex topological vector spaces. Then, in section 4, we prove that every Housdorrf quasi-point-separable topological vector space has the fixed point property, which is indeed a proper extension of the Tychonoff's fixed point theorem.

From section 3, the constructions of quasi-point-separable topological vector spaces are similar to that of pseudonorm adjoint topological vector spaces introduced in [9]. As a matter of fact, the key ideas of the extension are significant important under the following senses:

1. In general, it may not be very difficult to check the conditions for a topological vector space to be quasi-point-separable;
2. Every point-separable topological vector space is quasi-point-separable;
3. From the proof of Theorem 4.1 in section 4, the quasi-point-separable property is a sufficient condition for Housdorrf topological vector spaces to have the fixed point property;
4. The class of quasi-point-separable topological vector spaces is very large that includes pseudonorm adjoint topological vector spaces and locally convex topological vector spaces as special cases.

2. **Preliminaries**

As we mentioned in the introduction section, the concept of quasi-point-separable topological vector spaces is a generalization of pseudonorm adjoint topological vector spaces, which is a proper extension of locally convex topological vector spaces.  And the above fixed point theorem proved in [9] is a proper extension of the Tychonoff's fixed point theorem. For easy referee, in this section, we recall the concepts of pseudonorm adjoint topological vector spaces and its properties studied in [9].

**Definition 2.1** [9] Let $X$ be a vector space with origin $\theta$. A mapping $p: X \to \mathbb{R}^+$ is called a pseudonorm on $X$ if it satisfies the following conditions:

W₁. $p(x) \geq 0$, for all $x \in X$ and $p(\theta) = 0$;
W₂. $p(-x) = p(x)$, for all $x \in X$;
W₃. For any elements $x_1$, $x_2$ of $X$, and $0 \leq \alpha \leq 1$, one has

$$p(\alpha x_1 + (1 - \alpha)x_2) \leq \alpha p(x_1) + (1 - \alpha)p(x_2).$$

**Definition 2.2 [9].** Let $X$ be a vector space. A mapping $q: X \to \mathbb{R}^+$ is called a quasi-pseudonorm on $X$ if there are a pseudonorm $p$ on $X$ and a strictly increasing continuous function $\varphi: \mathbb{R}^+ \to \mathbb{R}^+$ such that

W₄. $q(x) \leq \varphi(p(x))$, for all $x \in X$;
W₅. $\varphi(0) = 0$.

Here, $q$ is said to be adjoint with the pseudonorm $p$ and the weighted function $\varphi$.

**Definition 2.3 [9].** Let $(X, \tau)$ be a topological vector space. If $X$ is equipped with a family of $\tau$-continuous quasi-pseudonorms $\{q_\lambda\}_{\lambda \in \Lambda}$ associated with a family of $\tau$-continuous pseudonorms $\{p_\lambda\}_{\lambda \in \Lambda}$ and a family of weighted functions $\{\varphi_\lambda\}_{\lambda \in \Lambda}$, then $(X, \tau)$ is called a pseudonorm adjoint topological vector space.

**Definition 2.4 [9].** A family of quasi-pseudonorms $\{q_\lambda\}_{\lambda \in \Lambda}$ equipped on a topological vector space $(X, \tau)$ is said to be total whenever, for $x \in X$, $q_\lambda(x) = 0$ holds, for every $\lambda \in \Lambda$, then it is necessary to have $x = \theta$. A pseudonorm adjoint topological vector space is said to be total if it is equipped with a total family of quasi-pseudonorms.

**Lemma 2.5 [9].** *Every locally convex topological vector space is a pseudonorm adjoint topological vector space.*

### 3. M-quasiconvex functions and quasi-point-separable topological vector spaces

**Definition 3.1.** Let $X$ be a vector space with origin $\theta$. A quasiconvex function $u: X \to \mathbb{R}^+$ is said to be *m-quasiconvex* if it satisfies $u(\theta) = 0$, and, for any elements $x_1$, $x_2$ of $X$, and $0 \leq \alpha \leq 1$, one has

$$u(\alpha x_1 + (1 - \alpha)x_2) \leq \text{Max}\{u(x_1), u(x_2)\}.$$

We list some useful examples of m-quasiconvex functions on $X$ below.

(a) *Every nonnegative convex function satisfying $u(\theta) = 0$ is m-quasiconvex;*
(b) *Every semi-norm is m-quasiconvex.*

**Lemma 3.3.** *Suppose that $u$ is an m-quasiconvex function on $X$. Then, for any $x \in X$, the function $u(tx)$ is an increasing function with respect to $t > 0$.*

*Proof.* Take arbitrary $0 < t_1 < t_2$. Then

$$u(t_1 x) = u\left(\left(1 - \frac{t_1}{t_2}\right)\theta + \frac{t_1}{t_2}t_2 x\right) \leq \text{Max}\{u(\theta), u(t_2 x)\} = u(t_2 x). \qquad \square$$

**Lemma 3.4.** *Let $(X, \tau)$ be a topological vector space with dual space $X^*$. Then, for every $u \in X^*$, $|u|$ is a $\tau$-continuous m-quasiconvex function on $X$.*

**Lemma 3.5.** *Every $\tau$-continuous pseudonorm $p: X \to \mathbb{R}^+$ is a $\tau$-continuous m-quasiconvex function on $X$.*

**Definition 3.6**. Let $(X, \tau)$ be a topological vector space with dual space $X^*$ (the space of linear and $\tau$-continuous functions on $X$).

(i) If there is a subset $V \subseteq X^*$, such that, for $x \in X$,

$$v(x) = 0, \text{ for every } v \in V, \text{ implies } x = \theta,$$

then $(X, \tau)$ is said to be point-separable (by $V$) (or the set $V$ is said to be total on $X$). Meanwhile, the set $V$ is called a point-separating space for this vector space $X$.

(ii) Suppose that a topological vector space $(X, \tau)$ is equipped with a family of $\tau$-continuous m-quasiconvex functions $\{u_\lambda\}_{\lambda \in \Lambda}$, where $\Lambda$ is an index set. If, for $x \in X$,

$$u(x) = 0, \text{ for every } \lambda \in \Lambda, \text{ implies } x = \theta,$$

then $(X, \tau)$ is said to be quasi-point-separable (or the family $\{u_\lambda\}_{\lambda \in \Lambda}$ is said to be total on $X$). Or we say that $X$ is quasi-point-separated by the family $\{u_\lambda\}_{\lambda \in \Lambda}$. Meanwhile, the family $\{u\}_{\lambda \in \Lambda}$ is called a quasi-point-separating space for $X$.

**Lemma 3.7**. *Let $(X, \tau)$ be a topological vector space with dual space $X^*$. If $X$ is point-separable (by a point-separable space $V \subseteq X^*$), then $X$ is quasi-point-separable with the quasi-point-separating space $\{|v|: v \in V\}$.*

**Lemma 3.8**. *Let $(X, \tau)$ be a Hausdorff locally convex topological vector space with dual space $X^*$. Then $X$ is point-separable (by $X^*$), so is quasi-point-separable with the quasi-point-separating space $\{|v|: v \in X^*\}$.*

*Proof.* Suppose that $(X, \tau)$ is equipped with a total family of $\tau$-continuous seminorms $\{p_\lambda\}_{\lambda \in \Lambda}$ such that the initial topology $\tau$ on $X$ is induced by $\{p_\lambda\}_{\lambda \in \Lambda}$. Notice that the totality of the family of seminorms $\{p_\lambda\}_{\lambda \in \Lambda}$ equipped on $(X, \tau)$ is equivalent to the Hausdorff property of $(X, \tau)$. □

**Lemma 3.9**. *Suppose that a total pseudonorm adjoint topological vector space $(X, \tau)$ is equipped with a total family of $\tau$-continuous quasi-pseudonorms $\{q_\lambda\}_{\lambda \in \Lambda}$ associated with a family of $\tau$-continuous pseudonorms $\{p_\lambda\}_{\lambda \in \Lambda}$ and a family of weighted functions $\{\varphi_\lambda\}_{\lambda \in \Lambda}$. Then $X$ is quasi-point-separable with the quasi-point-separating space $\{p_\lambda\}_{\lambda \in \Lambda}$.*

## 4. The fixed point property of quasi-point-separable Housdorrf topological vector spaces

We prove the main theorem of this paper in this section. The ideas of the proof of this theorem is similar to the proof of the main theorem in [9]. As mentioned in section 2 and as in the proof of this theorem, we see that the quasi-point-separable property of topological vector spaces plays the key factor for the considered spaces to have the fixed point property.

**Theorem 4.1.** *Every quasi-point-separable Housdorrf topological vector space has the fixed point property.*

*Proof.* Let $(X, \tau)$ be a quasi-point-separable Housdorrf topological vector space with a quasi-point-separating space $\{u_\lambda\}_{\lambda \in \Lambda}$. Let $C$ be a nonempty compact convex subset of $X$. Let $f: C \to C$ be a continuous mapping. The point-separating property of $\{u_\lambda\}_{\lambda \in \Lambda}$ implies that a point $x_0 \in C$ is a fixed point of $f$ if and only if

$$u_\lambda(x_0 - f(x_0)) = 0, \text{ for every } \lambda \in \Lambda.$$

Therefore, it is sufficient to show

$$\bigcap_{\lambda \in \Lambda}\{x \in C: u_\lambda(x - f(x)) = 0\} \neq \emptyset. \tag{4.1}$$

To this end, we first show that, for an arbitrary positive integer $m$ and an arbitrary finite subset $\{\lambda_1, \lambda_2, \ldots, \lambda_m\} \subseteq \Lambda$, the following statement holds:

$$\bigcap_{1 \leq k \leq m}\{x \in C: u_{\lambda_k}(x - f(x)) = 0\} \neq \emptyset. \tag{4.2}$$

Assume contrarily that (4.2) does not hold. That is,

$$\bigcap_{1 \leq k \leq m}\{x \in C: |u_{\lambda_k}(x - f(x))| = 0\} = \emptyset. \tag{4.3}$$

Take a strictly decreasing sequence of positive numbers $\{\delta_i\}$ with limit 0; i.e., $\delta_i \downarrow 0$, as $i \to \infty$. For every $i$, let

$$G_i = \bigcap_{1 \leq k \leq m}\{x \in C: |u_{\lambda_k}(x - f(x))| \leq \delta_i\}$$
$$= \{x \in C: \text{Max}_{1 \leq k \leq m} |u_{\lambda_k}(x - f(x))| \leq \delta_i\}.$$

The $\tau$-continuity of $f: C \to C$ and the $\tau$-continuity of the m-quasiconvex functions $u_{\lambda_k}$ imply that $\{G_i\}$ is a decreasing (with respect to inclusion order) sequence of $\tau$-closed subsets of the compact set $C$. We have

$$\bigcap_{1 \leq k \leq m}\{x \in C: |u_{\lambda_k}(x - f(x))| = 0\} = \bigcap_{1 \leq i < \infty} G_i$$

By the hypothesis (4.3), there is a positive integer $N$ such that $G_i = \emptyset$, for all $i \geq N$. It follows that there is $\delta > 0$ such that

$$\text{Max}_{1 \leq k \leq m} |u_{\lambda_k}(x - f(x))| \geq \delta, \text{ for every } x \in C. \tag{4.4}$$

Based on the inequality (4.4), we define a set-valued mapping $F: C \to 2^C \setminus \{\emptyset\}$ as follows:

$$F(x) = \{z \in C: \text{Max}_{1 \leq k \leq m} |u_{\lambda_k}(x - f(z))| \geq \delta\}, \text{ for } x \in C.$$

Then, for every $x \in C$, we have

$$F(x) = \bigcup_{1 \leq k \leq m}\{z \in C: |u_{\lambda_k}(x - f(z))| \geq \delta\}.$$

For every $x \in C$, from (4.4), $x \in F(x)$, and from the $\tau$-continuity of $f: C \to C$ and the $\tau$-continuity of the m-quasiconvex functions $u_{\lambda_k}$ again, it implies that, $F(x)$ is a nonempty $\tau$-closed subset of $C$.

Next, we show that the mapping $F: C \to 2^C \setminus \{\emptyset\}$ is a KKM mapping. For any given positive integer $n$, take arbitrary $n$ points $x_1, x_2, \ldots, x_n \in C$. For any positive numbers $t_1, t_2, \ldots, t_n$ satisfying $\sum_{i=1}^n t_i = 1$, let $y = \sum_{i=1}^n t_i x_i$. We show that

$$y \in \bigcup_{1 \leq i \leq n} F(x_i). \tag{4.5}$$

Assume, by the way of contradiction, that (4.5) does not hold. Then we have

$$y \notin F(x_i), \text{ for every } i = 1, 2, \ldots, n.$$

It is

$$\text{Max}_{1 \leq k \leq m} |u_{\lambda_k}(x_i - f(y))| < \delta, \text{ for all } i = 1, 2, \ldots, n. \tag{4.6}$$

From assumption (4.4), the inequality (4.6), and by the definition of m-quasiconvex functions $u_{\lambda_k}$, it follows that

$$\delta \leq \text{Max}_{1 \leq k \leq m} |u_{\lambda_k}(y - f(y))|$$

$$= \text{Max}_{1 \leq k \leq m} |u_{\lambda_k}(\sum_{i=1}^n t_i x_i - f(y))|$$

$$= \text{Max}_{1 \leq k \leq m} |u_{\lambda_k}(\sum_{i=1}^n t_i (x_i - f(y)))|$$

$$\leq \text{Max}_{1 \leq k \leq m} (\text{Max}_{1 \leq i \leq n} |u_{\lambda_k}(x_i - f(y))|)$$

$$= \text{Max}_{1 \leq i \leq n} (\text{Max}_{1 \leq k \leq m} |u_{\lambda_k}(x_i - f(y))|)$$

$$< \text{Max}_{1 \leq i \leq n} \delta$$

$$= \delta.$$

It is a contradiction. It implies that $F: C \to 2^C \setminus \{\emptyset\}$ is a KKM mapping with nonempty $\tau$-closed values. Since $C$ is compact, from Fan-KKM Theorem, we have

$$\bigcap_{x \in C} \{z \in C: \text{Max}_{1 \leq k \leq m} |u_{\lambda_k}(x - f(z))| \geq \delta\} = \bigcap_{x \in C} F(x) \neq \emptyset.$$

Then, there is $z_0 \in \bigcap_{x \in C} F(x)$ satisfying

$$\text{Max}_{1 \leq k \leq m} |u_{\lambda_k}(x - f(z_0))| \geq \delta, \text{ for every } x \in C.$$

In particularly, if we take $x = f(z_0)) \in C$, we get

$$0 = \text{Max}_{1 \leq k \leq m} |u_{\lambda_k}(f(z_0) - f(z_0))| \geq \delta.$$

So it is a contradiction to the hypothesis (4.3). It implies that for any positive integer $m$ and arbitrary finite subset of $\Lambda$ $\{\lambda_1, \lambda_2, \ldots, \lambda_m\}$, we must have

$$\bigcap_{1 \leq k \leq m} \{x \in C : |u_{\lambda_k}(x - f(x))| = 0\} \neq \emptyset.$$

Then (4.1) immediately follows from the finite nonempty intersection property and the compactness of $C$. □

In [13], Park studied fixed point problems on point-separable topological vector spaces for some classes of functions, such as half-continuous functions. Here, we use Theorem 4.1 to obtain the following fixed point theorem. Since every Hausdorff locally convex topological vector space is a point-separable topological vector space, the following theorem is also an extension of the Tychonoff's fixed point theorem to quasi-point-separable topological vector spaces.

**Theorem 4.3**. *Every point-separable topological vector space has the fixed point property.*

**Theorem 4.1 [9] (an alternative version).** *Every Housdorrf and total pseudonorm adjoint topological vector space has the fixed point property.*

**Tychonoff's fixed point theorem** [18] (**an alternative version**). *Every Hausdorff locally convex topological vector space has the fixed point property.*

## 5. Examples

In 1935, Tychonoff proved that the topological vector space $l_r$, for $0 < r < 1$, is not locally convex (see [4], [7−10], [18]). In [9], it is proved that the topological vector space $l_r$, for $0 < r < 1$, is a Hausdorff pseudonorm adjoint topological vector space.

In this section, one-step further, we give a simple proof showing that, for $0 < r < 1$, the metric vector spaces $l_r$ and $l^r$ (Examples 5.1 and 5.2, respectively) are point-separable; and therefore, from Lemma 3.7 or 3.8, it follows immediately that $l_r$ and $l^r$ both are quasi-point-separable, which are not locally convex.

Then, from Lemma 3.4 in [9] and Lemma 3.8, every Hausdorff locally convex topological vector space is a Hausdorff quasi-point-separable topological vector space. Hence, Theorem 4.1 also properly extends Tychonoff fixed point theorem.

In this section, let $S$ denote the set of sequences of real numbers.

**Example 5. 1**. For every given $r \in (0, 1)$, define a subspace $l_r$ of $S$ as below

$$l_r = \{\{x_i\} \in S : \sum_{i=1}^{\infty} |x_i|^r < \infty\}.$$

Define a function $p$ on $l_r \times l_r$ as

$$p(\{x_i\}, \{y_i\}) = \sum_{i=1}^{\infty} |x_i - y_i|^r, \text{ for } \{x_i\}, \{y_i\} \in l_r. \tag{5.1}$$

Then $(l_r, p)$ is a metric vector space (it is not locally convex) and it has the fixed point property.

*Proof.* For every fixed $n = 1, 2, \ldots$, define a function $v_n$ on $l_r$ as follows:

$$v_n(\{x_i\}) = x_n, \text{ for } \{x_i\} \in l_r. \tag{5.2}$$

Then $\{v_n\} \subseteq l_r^*$. We see that, for $\{x_i\} \in l_r$,

$$v_n(\{x_i\}) = 0, \text{ for all } n = 1, 2, \ldots, \quad \text{implies} \quad \{x_i\} = \theta.$$

So the space $(l_r, p)$ is point-separated by $\{v_n\} \subseteq l_r^*$; and therefore, from Lemma 3.8, this metric space $(l_r, p)$ is quasi-point-separable. By Theorem 4.1, this space $(l_r, q)$ has the fixed point property. □

**Example 5. 2.** For every given $r \in (0, 1)$, define a subspace $l^r$ of $S$ as below

$$l^r = \left\{\{x_i\} \in S : \sum_{i=1}^{\infty} \frac{|x_i|^r}{1+|x_i|^r} < \infty \right\}.$$

Define a function $q$ on $l^r \times l^r$ as

$$q(\{x_i\} - \{y_i\}) = \sum_{i=1}^{\infty} \frac{|x_i - y_i|^r}{1+|x_i - y_i|^r}, \text{ for every } \{x_i\}, \{y_i\} \in l^r. \tag{5.3}$$

Then $(l^r, q)$ is a metric vector space (it is not locally convex) and it has the fixed point property.

*Proof.* To show that the functional $q$ defined in (5.3) defines a metric on $l^r$, one may consider the following steps:

(a) $\frac{t^r}{1+t^r}$ is a strictly increasing function of $t$ on $\mathbb{R}^+$;
(b) $(t+s)^r \leq t^r + s^r$, for $t, s \in \mathbb{R}^+$ (this should be useful in Example 5.1);
(c) $\frac{(t+s)^r}{1+(t+s)^r} \leq \frac{t^r}{1+t^r} + \frac{s^r}{1+s^r}$, for $t, s \in \mathbb{R}^+$.

For every fixed $n = 1, 2, \ldots$, similarly to Example 5.1, define a function $w_n$ on $l^r$ as follows:

$$w_n(\{x_i\}) = x_n, \text{ for } \{x_i\} \in l^r. \tag{5.2}$$

Then, $\{w_n\} \subseteq l^{r*}$. We see that, for $\{x_i\} \in l^r$,

$$w_n(\{x_i\}) = 0, \text{ for all } n = 1, 2, \ldots, \quad \text{implies} \quad \{x_i\} = \theta.$$

Rest of the proof will be similar to the proof of Example 5.1. □